\documentclass[12pt]{amsart}
\usepackage{verbatim}

\usepackage[english]{babel}
\usepackage{inputenc}
\usepackage{amssymb,amsmath,amsthm,amscd, enumerate,mathrsfs,graphicx,colordvi}
\usepackage[all]{xy}
\usepackage{color}   
\newcommand{\rc}{\color{red}} 
\newcommand{\ec}{\color{black}}
\newcommand{\bc}{\color{black}}
\newcommand{\gc}{\color{green}}

\newcommand{\f}{\mathfrak}
\newcommand{\ms}{\mathscr}
\newcommand{\mc}{\mathcal}
\newcommand{\ff}{\ms F}
\newcommand{\uu}{\ms U}
\newcommand{\nilp}{{\rm Nilp}}
\newcommand{\q}{{J}}
\newcommand{\p}{{P}}
\newcommand{\n}{{\mathfrak Q}}
\newcommand{\s}{{\mathcal E}}
\newcommand{\m}{{M}}
\newcommand{\aaa}{{I}}
\newcommand{\nn}{{N}}
\newcommand{\ii}{{\cdot}}
\newcommand{\pp}{{\mathfrak P}}
\newcommand{\fa}{{A\Join^n\!\q}}
\newcommand{\da}{{A\!\Join^f\!\!\f b}}
\newcommand{\dda}{{D\!\Join^f\!\!J}}
\newcommand{\lam}{\lambda}
\newcommand{\La}{\Lambda}
\newcommand{\spec}{{\rm Spec}}
\newcommand{\zar}{{\rm Zar}}
\newcommand{\ad}{{\rm Ad}}
\newcommand{\jac}{{\rm Jac}}%
\newcommand{\rad}{{\rm rad}}%
\newcommand{\kr}{{\rm Kr}}
\newcommand{\hgt}{{\rm ht}}
\newcommand{\Max}{{\rm Max}}
\newcommand{\clop}{{\rm Clop}}
\newcommand{\stf}{\star{_{_{\!}f}}}
\newcommand{\Min}{{\rm Min}}
\newcommand{\Ima}{{\rm Im}}
\newcommand{\Ker}{{\rm Ker}}
\newcommand{\tot}{{\rm Tot}}
\newcommand{\dd}{{\displaystyle}}
\newcommand{\z}{{\ldots}}
\newcommand{\w}{{\setminus}}
\newcommand{\e}{\grave{e}}
\newcommand{\ude}{\mbox{\textsl{id}}} %
\newcommand{\udes}{\mbox{\textsl{\tiny id}}} %

\theoremstyle{break}
\newtheorem{thm}{ \textbf{Theorem}}[section]
\newtheorem{defn}[thm]{ \textbf{Definition}}
\newtheorem{cor}[thm]{ \textbf{Corollary}}
\newtheorem{lem}[thm]{ \textbf{Lemma}}
\newtheorem{prop}[thm]{ \textbf{Proposition}}

\theoremstyle{definition}
\newtheorem{ex}[thm]{ \textbf{Example}}
\newtheorem{oss}[thm]{Remark}
\theoremstyle{remark}

\title[Pr\"ufer-like conditions]
{Pr\"ufer--like conditions on an amalgamated algebra along an ideal}
\author[C. A. Finocchiaro]
{Carmelo Antonio Finocchiaro}
\address{Dipartimento di Matematica -- Universit\`{a} degli Studi Roma Tre, Largo San Leonardo Murialdo 1, 00146 Roma}



\begin{document}

\maketitle

\begin{abstract}
Let $f:A\longrightarrow B$ be a ring homomorphism and let $\f b$ be an ideal
 of $B$. In this paper we study Pr\"ufer--like conditions in the amalgamation of $A$
  with $B$ along $\f b$, with respect to $f$, a ring construction introduced in 2009 by D'Anna, Finocchiaro and Fontana.
\end{abstract}
 
\begin{flushright}
\scriptsize \emph{In memory of my father}
\end{flushright}
\section{Introduction}
 Pr\"ufer domains, introduced by H. Pr\"ufer in \cite{Pr}, form a very relevant class of  commutative rings. Throughout the years, this class was deeply studied by several authors (for a sistematic study see \cite{fo-hu-pa}), so that many equivalent definitions of a Pr\"ufer domain were given. For example, the notion of Pr\"ufer domain globalizes the notion of valuation domain in a non local context. Moreover, the class of Pr\"ufer domains is the natural generalization of the class of Dedekind domains in the non--Noetherian setting.
Among the many equivalent conditions that make an integral  domain $A$ a Pr\"ufer domain, we recall the following:
\begin{enumerate}
\item[(1)] Every finitely generated ideal of $A$ is projective.
\item[(2)] $A_{\f p}$ is a valuation domain, for each prime (maximal) ideal $\f p$ of $A$.
\item[(3)] Every finitely generated ideal of $A$ is locally principal.
\item[(4)] If $T$ is an indeterminate over $A$, every polynomial $f\in A[T]$ is a  Gauss polynomial over $A$ (i.e., $c(fg)=c(f)c(g)$, for each polynomial $g\in A[T]$, where $c(f)$ denotes the content of the polynomial $f$).
\item[(5)] Every nonzero finitely generated ideal of $A$ is invertible.
\end{enumerate}
In \cite{Gr}, the notion of  Pr\"ufer domain was generalized to arbitrary (commutative) rings, possibly with zerodivisors.
Other important contributions to the study of the previous conditions in rings with zerodivisors were
given in \cite{B-Gl-2006}, \cite{endo}, \cite{fu}, \cite{Glaz-2005-1}, \cite{Glaz-2005-2}, \cite{Je}, \cite{Lu-Gaussian},
\cite{Lu}, \cite{Tsang}, etc. On the other hand, recently, in \cite{Ba-Gl}, Bazzoni and Glaz showed,
giving appropriate counterexamples, that none of the previous conditions is equivalent to the other, when $A$ is a ring with zerodivisors.
The fact that, in general, the rings satisfying previous Pr\"ufer--like conditions are distinct classes
of rings leads us to recall the following definition.
\begin{defn}
Let $A$ be a ring.
\begin{enumerate}[$(P_1)$]
\item $A$ is a \emph{semi--hereditary ring} if every finitely generated ideal of $A$ is projective.
\item $A$ \emph{has weak global dimension at most 1} if $A_{\f p}$ is a valuation domain, for each prime (maximal) ideal $\f p$ of $A$.
\item $A$ is an \emph{arithmetical ring} if every finitely generated ideal of $A$ is locally principal.
\item $A$ is a \emph{Gauss ring} if every polynomial $f\in A[T]$ is a  Gauss polynomial over $A$.
\item $A$ is a \emph{Pr\"ufer ring} if every regular and finitely generated ideal of $A$ is invertible.
\end{enumerate}
\end{defn}
In \cite{Ba-Gl}, it is shown that, for each $n\in\{1,2,3,4\}$, condition (P$_n$) implies condition (P$_{n+1}$). More precisely, Bazzoni and Glaz proved that a ring $A$ satisfies condition (P$_n$) if and only if $A$ satisfies condition (P$_{n+1})$ and the total ring of fractions $\tot(A)$ of $A$ satisfies condition (P$_n$). Moreover, it is proved that, if $\tot(A)$  is an  absolutely flat ring, then conditions (P$_n$) ($n\in\{1,2,3,4,5\}$) are equivalent on $A$.

Recently, J. Boynton in \cite{boy} studied Pr\"ufer--like conditions in pullbacks. The use of pullbacks and fiber products of ring homomorphisms is a very powerful tool to produce interesting examples (see \cite{fo}, \cite{ga}, \cite{ga-ho}). Of particular interest are the pullbacks of the following type: let $A\subseteq B$ be a ring extension such that $A$ and $B$ have a nonzero common ideal. In this case, call \emph{the conductor of $B$ into $A$} the largest nonzero common ideal to $A$ and $B$. It is well-known that the conductor of such a ring extension $A\subseteq B$ is
$$
\f c:=(B:A):=\{x\in A:xB\subseteq A\}
$$
Thus, if $\pi:B\longrightarrow B/\f c$ is the canonical projection, then $A$ is clearly the inverse image $\pi^{-1}(A/\f c)$ of the subring $A/\f c$ of $B/\f c$.

In his paper, Boynton describes the transfer of Pr\"ufer--like conditions on this kind of pullbacks, under the assumption that the conductor of $A\subseteq B$ is a regular ideal of $B$.

The aim of the present paper is to study Pr\"ufer--like conditions on amalgamated algebras along ideals. More precisely, in \cite{da-fi-fo1} and \cite{da-fi-fo2}, the authors have introduced  the following new ring construction. Given a ring homomorphism $f:A\longrightarrow B$ and an ideal $\f b$ of $B$, consider the subring
$$
\da:=\{(a,f(a)+b):a\in A,b\in\f b\}
$$
of $A\times B$, called \emph{the amalgamation of $A$ with $B$ along $\f b$ with respect to $f$}. This construction generalizes the amalgamated duplication of a ring along an ideal (introduced and studied in \cite{d'a}, \cite{d'a-f-1},
\cite{d'a-f-2} and in \cite{m-y}). Moreover, several classical constructions (such as the $A+XB[X]$, the
$A+XB[\![X]\!]$ and the $D+M$ constructions) can be studied as particular cases of the amalgamation (see \cite[Examples 2.5 and 2.6]{da-fi-fo1}) and other classical constructions, such as the Nagata's idealization (cf. \cite[page
2]{N}, \cite[Chapter VI, Section 25]{h}), and the CPI extensions (in the sense of Boisen and Sheldon \cite{bo}) are related to it (see \cite[Example 2.7]{da-fi-fo1}).
The level of generality choosen to define the amalgamation is due to the fact that the ring $\da$ may be studied in the frame of fiber product constructions. This allows to describe easily many algebraic properties of $\da$, in relation with those of $A,B,\f b$ and $f$.

Moreover, the ring $A$ is always embedded into the ring $\da$, and the natural image of the ring $A$ into $\da$ is a retract of $\da$ (see \cite[Remark 4.6 or Proposition 4.7]{da-fi-fo1}. This will help us to describe the transfer of Pr\"ufer--like conditions in the amalgamations.

The content of this paper is organized as follows: at the beginning, we prove that, under the assumption that the conductor of the ring extension $\da \subseteq A\times B$ is regular, the ring  $\da$ satisfies Pr\"ufer--like conditions (P$_n$), for $n\in\{1,2,3,4,5\}$, only in the trivial cases.

Later, we investigate the general case (in which the conductor is not necessarily regular) and we provide sufficient and necessary conditions for $\da$ to satisfy conditions (P$_n$), for $n\in\{1,2,3,4,5\}$.

The results of this paper form a part of author's thesis. The author is grateful to Marco Fontana, Stefania Gabelli and Sarah Glaz for their helpful comments and suggestions.

\section{Preliminaries}
We begin with some terminology and notation.
In the following, with the term \emph{ring} we will mean a commutative ring with multiplicative identity. We will call an element of a ring $A$ a \emph{regular} element if it is not a zerodivisor, and set
$$
{\rm Reg}(A):=\{a\in A: a \mbox{ is a regular element of }A\}
$$
Moreover, we will say that an ideal of $A$ is a \emph{regular ideal} if it contains a regular element of $A$. As usual, we will denote be $\spec(A)$ the set of all prime ideals of $A$ and sometimes, but not always, it will be endowed with the Zariski topology.

We collect in the following proposition several properties of the ring construction $\da$, that follow easily from the definitions.
\begin{prop}{\rm(\cite[Proposition 5.1]{da-fi-fo1})}
\label{inizio}   Let $f: A\longrightarrow B$ be a ring   homomorphism,  $\f b$ an ideal of $B$  and let
$$
\da := \{ (a,
f(a)+b) : a\in A, \ b \in \f b \}
$$
The following statements hold.
\begin{enumerate}
 \item[\rm (1)]  Let  $\iota := \iota_{A, f, \f b}: A\longrightarrow \da$ be the natural
 ring homomorphism defined by $\iota(a)  := (a, f(a))$, for all $a \in A$. Then, $\iota$ is ring embedding, making
$\da$   a  ring   extension of $A$ \  (with $\iota(A) =  \Gamma(f) \  (:=\{(a,f(a)) : a\in A\}$  subring of
$\da$).

\item[\rm (2)]  Let $\f a$ be an ideal of $A$ and set $ \f a\!\Join^f\!\! \f b :=\{(a, f(a)+b) : a\in \f a, b \in \f b \}$.
Then $\f a\!\Join^f\!\! \f b$  is an ideal of $\da$, the composition  of canonical homomorphisms
$A\stackrel{\iota}{\longrightarrow} \da\longrightarrow \da/\f a\!\Join^f\!\! \f b$ is a surjective ring homomorphism
and its kernel coincides with  $\f a$.\\ Hence, we have the following canonical isomorphism:
 $$
\frac{\da}{\f a\!\Join^f\!\! \f b}\cong \frac{A}{\f a}\,.
$$

 \item[\rm (3)] Let $p_{_A}: \da \longrightarrow A$ and $p_{_B}:\da\longrightarrow B$ be the natural projections of
 $\da \subseteq A \times B$ {into} $A$ and $B$, respectively.  Then, $p_{_A}$ is surjective
 and\, $\Ker(p_{_A})=\{0\}\times \f b$.\\
Moreover, $p_{_B}(\da)=f(A)+ \f b$ and $\Ker(p_{_B})=f^{-1}(\f b)\times \{0\}$.  Hence, the following canonical
isomorphisms  hold:
$$
\frac{\da}{(\{0\}\times \f b)}\cong A\quad \mbox{ and } \quad\frac{\da}{f^{-1}(\f b)\times \{0\}}\cong f(A)+\f b \,.
$$
\item[\rm(4)] Let $\gamma:\da\longrightarrow (f(A)+\f b)/\f b$
be the natural ring homomorphism, defined by $(a,f(a)+b)\mapsto f(a)+\f b$. Then $\gamma$ is surjective and
$\Ker(\gamma)=f^{-1}(\f b)\times \f b$. Thus, we have the following natural isomorphisms
$$
\frac{\da}{f^{-1}(\f b)\times \f b}\cong\frac{f(A)+\f b}{\f b}\cong \frac{A}{f^{-1}(\f b)}\,.
$$
In particular, when $f$ is surjective we have
$$
\frac{\da}{f^{-1}(\f b)\times \f b}\cong \frac{B}{\f b}\,.
$$
 \end{enumerate}
\end{prop}
\begin{defn}\label{fibproddef}
Let $\rho:A\longrightarrow C,\sigma:B\longrightarrow C$ be ring homomorphisms. We recall that the following subring
$$
\rho\times_C\sigma:=\{(a,b)\in A\times  B:\rho(a)=\sigma(b)\}
$$
of $A\times B$ is usually called \emph{the fiber product of $\rho$ and $\sigma$}.
\end{defn}
\begin{prop}[{\rm\cite[Proposition 4.2]{da-fi-fo1}}]\label{fibprodfarf} Let $f:A\longrightarrow B$ be a ring homomorphism, $\f b$ be an ideal of $B$. If $\pi:B\longrightarrow B/\f b$ is the canonical projection and $\check{f}:=\pi\circ f$, then $\da=\check{f}\times_{B/\f b}\pi$.
\end{prop}
{
\begin{oss}\label{b-loc}
Let $f:A\longrightarrow B$ be a ring homomorphism, $S$ be a multiplicative subset of $A$ and $\f b$ be an ideal of $B$. Consider the multiplicative subset $T:=f(S)+\f b$ of $B$ and let $f_S:A_S\longrightarrow B_T$ be the ring homomorphism induced by $f$. By a straightforward verification it is shown that ${f_S}^{-1}(\f b B_T)=f^{-1}(\f b)A_S$. Moreover, for each ideal $\f d$ of $B$, it is immediate that $\f d B_T=B_T$ if and only if $f^{-1}(\f b+\f d)\cap S\neq \emptyset$. Thus, $B_T=\{0\}$ if and only if $f^{-1}(\f b)\cap S\neq \emptyset$.

If $\f p$ is a prime ideal of $A$ and $S:=A\w \f p, T:=S_{\f p}:=f(S)+\f  b$, we shall denote $f_S$ simply by $f_{\f p}$ and $\f b B_{T}$ simply by $\f b_{S_{\f p}}$.
\end{oss}
}
The following result describes completely  the prime spectrum of $\da$.
\begin{prop}{\rm (\cite[Proposition 2.6]{da-fi-fo2} and \cite[Propositions 4.1 and 4.2]{da-fi-fo3})}\label{spec}
We preserve the notation of Proposition \ref{inizio}. Set
$
X:=\spec(A),Y:=\spec(B),W:=\spec(\da).
$
For each prime ideal $\f p$ of $A$ and each prime ideal $\f q$ of $B$ not containing $\f b$, set
$$
\f p^{\prime_f}:=\{(p,f(p)+b):p\in \f p, b\in \f b\}
$$
$$
\overline{\f q}^f:=\{(a,f(a)+b):a\in A,b\in\f b, f(a)+b\in \f q\}.
$$
Then, the following statements hold.
\begin{enumerate}
\item [\rm (1)] The map $\f p \mapsto \f p^{\prime_{_{\! f}}}$ establishes a closed embedding of $X$ into  $W$,
 so its image, which coincides with $V(\f b_0)$, is homeomorphic to $X$.
\item [\rm (2)] The map $\f q \mapsto\overline{\f q}^{_{_{f}}} $ is a homeomorphism of $Y\w V(\f b)$ onto
$W\w V(\f b_0)$.

\item [\rm (3)] The prime ideals of $\da$ are of the type\ $\f p^{\prime_{_{\! f}}}$ or  $\overline{\f q}^{_{_{ f}}}$\!,\,  for
 $\f p$ varying in $ X$ and $\f q$ in $Y\w V(\f b)$.
\item [\rm(4)] Let $\f p \in\spec(A)$.  Then,  $\f p^{\prime_{_{\! f}}}$
 is a maximal ideal of $\da$ if and only if $\f p$ is a maximal ideal of $A$.
  \item [\rm(5)] Let $\f q$ be a prime ideal of $B$ not containing $\f b$. Then,
   $\overline{ \f q}^{_{_{f}}}$ is a maximal ideal of $\da$ if and only if $\f q$ is a maximal ideal of $B$.\\
In particular:
$$
\Max(\da)=\{\f p^{{\prime_{_{\! f}}}} :  \f p\in\Max(A)\}  \cup \{\overline{\f q}^{_{_{f}}} : \f q \in \Max(B)\w V(\f b)\}.
$$
\item[\rm(6)] The ring $\da$ is local if, and only if, $A$ is local and $\f b\subseteq \jac(B)$. In this case, if $\f m$ is the maximal ideal of $A$, then the maximal ideal of $\da$ is $\f m^{{\prime_{_{\! f}}}}$.
In particular, if $A$ and $B$ are local rings and $\f b$ is a proper ideal of $B$, then $\da$ is a local ring.
\item[\rm (7)]
\begin{enumerate}[\rm (a)]
\item If $\f q$ is a prime ideal of $B$ not containing $\f b$, then $(\da)_{\overline{\f q}^f}$ is isomorphic to $B_{\f q}$.
\item If $\f p$ is a prime ideal of $A$, consider the multiplicative subset $S_{\f p}:=f(A\w \f p)+\f b$, and let $f_{\f p}:A_{\f p}\longrightarrow B_{S_{\f p}}$ be the ring homomorphism induced by $f$. Then $(\da)_{\f p^{\prime_f}}$ is isomorphic to $A_{\f p}\Join^{f_{\f p}}\f b_{S_{\f p}}$. In particular, if $\f p\nsupseteq f^{-1}(\f b)$, we have $B_{S_{\f p}}=\{0\}$ and thus $(\da)_{\f p^{\prime_f}}$ is isomorphic to $A_{\f p}$.
\end{enumerate}
\end{enumerate}
\end{prop}
{
Let $A$ be a ring and $\f p$ be a prime ideal of $A$. Recall that $(A,\f p)$ has \emph{the regular total order property} if, for each pair of ideals $\f a_1,\f a_2$ of $A$, one at least of which is regular, the ideals $\f a_1A_{\f p},\f a_2A_{\f p}$ are comparable.
}

The following characterization of Pr\"ufer rings will be useful. We recall it here for the reader convenience.
\begin{thm}\label{one}
{\rm (\cite[Theorem 13]{Gr})}
Let $A$ be a ring. The following conditions are equivalent.
\begin{enumerate}[\rm (i)]
\item $A$ is a Pr\"ufer ring.
\item If $\f a,\f b,\f c$ are ideals of $A$ and $\f b$ or $\f c$ is regular, then
$$
\f a(\f b\cap \f c)=\f a\f b\cap \f a\f c.
$$
\item { For each maximal ideal $\f m$ of $A$, $(A,\f m)$ has the regular total order property. }
\end{enumerate}
\end{thm}
\begin{defn}
We say that a ring $A$ is {\rm a locally Pr\"ufer ring} if $A_{\f m}$ is a Pr\"ufer ring, for each $\f m\in \Max(A)$.
\end{defn}
\begin{oss}\label{Pru-Lu}
Let $A$ be a ring.
\begin{enumerate}[(a)]
\item  By \cite[Proposition 2.10]{Lu}, if $A$ is a locally Pr\"ufer ring, then $A$ is a Pr\"ufer ring.
\item If $A$ is Gauss ring, then so is $A_{\f m}$, for each maximal ideal $\f m$ of $A$ (each localization of a Gauss ring is still a Gauss ring). It follows that $A$ is a locally Pr\"ufer ring.
\item Note that an example of a Pr\"ufer and non locally Pr\"ufer ring is given in \cite[Example 2.11]{Lu}. Moreover, as observed in \cite[Example 3.8]{Ba-Gl}, if $K$ is a field and $T_1,T_2$ are indeterminates over $K$, then $K[T_1,T_2]/(T_1,T_2)^3$ is a local total ring of fractions (and thus a locally Pr\"ufer ring) that is not a Gauss ring. Thus we have the following proper inclusions of classes of rings
$$
\{\mbox{Semihereditary rings}\}\subsetneq \{{\rm w.gl.dim}\leq 1\}\subsetneq \{\mbox{Arithmetical rings}\}\subsetneq
$$
$$
\subsetneq\{\mbox{Gauss rings}\}\subsetneq\{\mbox{Locally Pr\"ufer rings}\}\subsetneq\{\mbox{Pr\"ufer rings}\}
$$
\end{enumerate}
\end{oss}
\begin{oss}\label{Ba} Let $\{A_1,\z,A_r\}$ be a nonempty and finite collection of rings and let $A:=\prod_{i=1}^rA_i$. As noted by Bakkari in a recent preprint, posted on arXiv, for each $n\in\{1,2,3,4,5\}$, $A$ satisfies Pr\"ufer--like condition (P$_n$) if and only if $A_i$ satisfies the same Pr\"ufer--like condition (P$_n$), for each $i\in\{1,\z,r\}$.
\end{oss}
\section{Results when the conductor of the ring extension $\da\subseteq A\times B$ is regular}

As noted in \cite[Lemma 1.50]{fi}, the conductor of the ring extension $\da \subseteq A\times B$ is $\f c:=f^{-1}(\f b)\times \f b$.
The following results show that when $\f c$ is a regular ideal of $A\times B$ (i.e., if $f^{-1}(\f b), \f b$ are regular ideals of $A,B$, respectively), then $\da$ satisfies Pr\"ufer--like conditions (P$_n$) ($n\in \{1,2,3,4,5\}$) only in the trivial case.
{
\begin{thm}\label{prufer-like-reg}
Let $f:A\longrightarrow B$ be a ring homomorphism and let $\f b$ be an ideal of $B$. If $f^{-1}(\f b)$ and $\f b$ are regular ideals, then the following conditions are equivalent.
\begin{enumerate}[\rm (i)]
\item $\da$ is a Pr\"ufer ring.
\item $A,B$ are Pr\"ufer rings and $\f b=B$.
\end{enumerate}
\end{thm}
\begin{proof}
(ii)$\Longrightarrow$(i). By (ii), $\da=A\times B$. Then, it sufficies to apply \cite[Proposition 3]{Gi-Hu}.

(i)$\Longrightarrow$(ii). Assume, by contradiction, that $\f b$ is a proper ideal of $B$, and pick a maximal ideal $\f m$ of $A$ containing $f^{-1}(\f b)$. Consider the multiplicative subset $S_{\f m}:=f(A\w \f m)+\f b$ of $B$. By Proposition \ref{spec}(7), the localization of $\da$ at the maximal ideal
$$
\f m^{\prime f}:=\{(m,f(m)+b):m\in\f m,b\in\f b\}
$$
is isomorphic to $C:=A_{\f m}\Join^{f_{\f m}}\f b_{S_{\f m}}$ ($f_{\f m}:A_{\f m}\longrightarrow B_{S_{\f m}}$ is the ring homomorphism induced by $f$). Now, pick regular elements $a_0\in f^{-1}(\f b),b_0\in \f b$. Then, in particular, $\f a_1:=(a_0,b_0)\da$ is a regular ideal of $\da$. Set $a^*:=a_0/1 \in A_{\f m}$, $b^*:= b_0/1\in B_T$. Obviously, $a^*,b^*$ are regular elements.  Since $\da$ is a Pr\"ufer ring, $(\da,\f m^{\prime f})$ has the regular total order property, by Theorem \ref{one}. Thus, if $\f a_2:=(a_0,0)\da$,  the ideals
$$
(a^*,b^*)C=\f a_1C,\qquad (a^*,0)= \f a_2 C
$$
are comparable.  Since, in particular, $b^*\neq 0$, we have $(a^*,b^*)C\nsubseteq (a^*,0)C$. It follows that $(a^*,0)C\subseteq (a^*,b^*)C$. Thus, there exist  elements $\alpha\in A_{\f m},\beta\in\f bB_{S_\f m}$ such that
$$
(a^*,0)=(\alpha,f_{\f m}(\alpha)+\beta)(a^*,b^*)
$$
Keeping in mind that $a^*$ is regular, it follows that $\alpha=1$. Then, $b^*(1+\beta)=0$, and thus $\beta=-1$, since $b^*$ is regular. This implies $\f bB_{S_{\f m}}=B_{S_{\f m}}$, and, by Remark \ref{b-loc}, $f^{-1}(\f b)\nsubseteq \f m$, a contradiction. Thus $\f b=B$ and, consequently,  $\da=A\times B$. Then, the remaining part of statement (ii) follows by \cite[Proposition 3]{Gi-Hu}.
\end{proof}
\begin{cor}\label{cor-prufer-like-reg}
We preserve the notation of Proposition \ref{inizio} and let $n\in\{1,2,3,4,5\}$. If $f^{-1}(\f b)$ and $\f b$ are regular ideals, then the following conditions are equivalent.
\begin{enumerate}[\rm (i)]
\item $\da$ satisfies Pr\"ufer--like condition {\rm (P$_n$)} (resp. $\da$ is locally Pr\"ufer).
\item $A,B$ satisfy Pr\"ufer--like condition {\rm (P$_n$)} (resp. $A,B$ are locally Pr\"ufer rings) and $\f b=B$.
\end{enumerate}
\end{cor}
\begin{proof}
(ii)$\Longrightarrow$ (i). By (ii), $\da=A\times B$. Then, condition (i) follows by using Remark \ref{Ba}, \cite[Proposition 3]{Gi-Hu} and definitions.

(i)$\Longrightarrow$(ii). If $\da$ is locally Pr\"ufer, then it is a Pr\"ufer ring, by Remark \ref{Pru-Lu}(a). Thus $\f b=B$, by Theorem \ref{prufer-like-reg}. Moreover, it is immediately seen that $A,B$ are locally Pr\"ufer rings. Now, let $n\in\{1,2,3,4,5\}$. If $\da$ satisfies Pr\"ufer--like condition P$_n$, then $\da$ is a Pr\"ufer ring. Thus the conclusion follows by Theorem \ref{prufer-like-reg} and Remark \ref{Ba}.
\end{proof}
\begin{cor}
Let $A$ be a ring and $\f a$ be a regular ideal of $A$. Consider \emph{the amalgamated duplication of $A$ along $\f a$}
$$
A\Join \f a:=\{(a,a+\alpha):a\in A,\alpha\in \f a\}
$$
(see \cite{d'a}, \cite{d'a-f-1}, \cite{d'a-f-2}) and let $n\in\{1,2,3,4,5\}$. Then, $A\Join \f a$ satisfies Pr\"ufer--like condition {\rm (P$_n$)} (resp. $A\Join \f a$ is a locally Pr\"ufer ring) if and only if $A$ satisfies Pr\"ufer--like condition {\rm (P$_n$)} (resp. $A$ is a locally Pr\"ufer ring) and $\f a=A$.
\end{cor}
\begin{proof}
Apply Corollary \ref{cor-prufer-like-reg}, keeping in mind \cite[Example 2.4]{da-fi-fo1}.
\end{proof}
}
\section{Results in the general case}
\begin{lem}\label{gauss-ret}
Let $r:B\longrightarrow A$ be a ring retraction, and $T$ be an  indeterminate over $B$. If $\sum_{i=0}^nb_iT^i$ is a Gauss polynomial over $B$, then $\sum_{i=0}^nr(b_i)T^i$ is a Gauss polynomial over $A$.
\end{lem}
\begin{proof}
It follows by the proof of \cite[Theorem 2.1(1)]{Ba-Ma-Mo}.
\end{proof}
\begin{prop}
We preserve the notation of Proposition \ref{inizio}. If $\da$ is a Pr\"ufer ring and $f({\rm Reg}(A))\subseteq {\rm Reg}(B)$, then $A$ is a Pr\"ufer ring.
\end{prop}
\begin{proof}
Let $T$ be an indeterminate over $A$ and $\f a:=(a_0,\z,a_n)$ be a regular and finitely generated ideal of $A$. Consider the polynomial $p(T):=\sum_{i=0}^na_iT^i\in A[T]$. Pick a regular element $a\in \f a$. Then, keeping in mind that $f({\rm Reg}(A))\subseteq {\rm Reg}(B)$, it is easily checked that $(a,f(a))$ is a regular element of the finitely  generated ideal $\f a^{\Join}:=((a_0,f(a_0)),\z, (a_n,f(a_n))$ of $\da$. Since $\da$ is a Pr\"ufer ring, it follows that $\f a^{\Join}$ is an invertible ideal of $\da$, and thus the polynomial $p_{\Join}(T):=\sum_{i=0}^n(a_i,f(a_i))
T^i\in \da[T]$, whose content is clearly $\f a^{\Join}$, is a Gauss polynomial over $\da$, by \cite{Tsang}. Let $p_{_A}:\da\longrightarrow A$ be the projection ($(a,f(a)+b)\mapsto a$). Then we have  $p(T)=\sum_{i=0}^np_{_A}
((a_i,f(a_i)))T^i$. Since $p_{_A}$ is a ring retraction (\cite[Remark 4.6]{da-fi-fo1}), it follows that $p(T)$ is a Gauss  polynomial over $A$, by Lemma \ref{gauss-ret}. Thus its content, that is exactly the regular ideal $\f a$, is invertible, by \cite[Theorem 6]{Lu-Gaussian}. This completes the proof.
\end{proof}
\begin{oss}
We preserve notation of Proposition \ref{inizio}. The fact that $\da$ is a Pr\"ufer ring does not imply, in general, that $A$ is a Pr\"ufer ring. For an example, see \cite[Example 2.3]{Ba-Ma-Mo}, keeping in mind \cite[Remark 2.8]{da-fi-fo1}.
\end{oss}
The following result is obtained by modifing the proof of \cite[Theorem 1]{Bo-La}.
\begin{prop}\label{im-reg}
Let $\phi:A\longrightarrow B$ be a surjective ring homomorphism. If $A$ is a Pr\"ufer ring and $\Ker(\phi)$ is a regular ideal of $A$, then $\f a(\f b\cap \f c)=\f a\f b\cap \f a\f c$, for all ideals $\f a,\f b,\f c$ of $B$. In particular, $B$ is a Pr\"ufer ring.
\end{prop}
\begin{proof}
Let $\f d:=\Ker(\phi)$ and let $\f a,   \f b, \f c$ be ideals of $B$. To prove the equality  $\f a(\f b\cap \f c)=\f a\f b\cap \f a\f c$, it sufficies to show that $\f a\f b\cap \f a \f c\subseteq \f a(\f b\cap \f c)$. If $\overline x\in \f a\f b\cap \f a \f c$, then there are elements $\overline a_i\in \f a,\overline b_i\in \f b, \overline\alpha_j\in\f a,\overline c_j \in \f c$, with $i\in\{1,\z,n\},j\in\{1,\z,m\}$, such that $\overline x=\sum_{i=1}^n\overline a_i\overline b_i=\sum_{j=1}^m\overline\alpha_j\overline c_j$. For each $i\in\{1,\z,n\}, j\in\{1,\z,m\}$, choose elements $a_i\in \phi^{-1}(\overline a_i),b_i\in \phi^{-1}(\overline b_i), \alpha_j\in \phi^{-1}(\overline \alpha_j), c_j\in \phi^{-1}(\overline c_j)$, and set $\f a':=\phi^{-1}(\f a),\f b':=\phi^{-1}(\f b),\f c':=\phi^{-1}(\f c)$. If $x:=\sum_{i=1}^na_ib_i$, it is immediate that $x-\sum_{j=1}^m\alpha_jc_j\in \f d$. Therefore $x\in (\f a'\f c'+\f d)\cap \f a'\f b'$. Keeping in mind Theorem \ref{one} and the fact that $\f d$ is a regular ideal of $A$, we have
$$
(\f a'\f c'+\f d)\cap \f a'\f b'=
(\f a'\f c' \cap \f a'\f b')+(\f d\cap \f a'\f b')\subseteq \left(\f a'\f c'\cap \f a'(\f b'+\f d)\right)+\f d=
$$
$$=\f a'(\f c'\cap(\f b'+\f d))+\f d
$$
Thus, there are elements $a'_h\in \f a',b'_h\in\f b',d_h\in\f d$, with $h\in\{1,\z,r\}$ such that $b'_h+d_h\in \f c'$, for each $h$, and $x=\sum_{h=1}^ra'_h(b'_h+d_h)+d$, for some $d\in \f d$. It follows immediately that $\overline x=\sum_{h=1}^rf(a'_h)f(b'_h)\in \f a(\f b\cap \f c)$. Now the first statement is clear. The fact that $B$ is a Pr\"ufer ring follows by the previous statement and Theorem \ref{one}.
\end{proof}
{\bc
\begin{cor}
Preserve the notation of Proposition \ref{inizio} and assume that  $\da$ is a Pr\"ufer ring. Then the following statements hold.
\begin{enumerate}[\rm (1)]
\item If $\{0\}\times \f b$ is a regular ideal of $\da$, then $A$ is a Pr\"ufer ring.
\item If $f^{-1}(\f b)\times\{0\}$ is a regular ideal of $\da$, then $f(A)+\f b$ is a Pr\"ufer ring.
\end{enumerate}
\end{cor}
\begin{proof}
It sufficies to apply Propositions \ref{inizio}(3) and \ref{im-reg}.
\end{proof}
Now, we will give sufficient conditions to make $\da$ a total ring of fractions (and, in particular, a Pr\"ufer ring).
\begin{prop}\label{Prufer-sufficient}
Let $A$ be a total ring of fractions (i.e. $A=\tot(A)$), $f:A\longrightarrow B$ be a ring homomorphism and $\f b$ be an ideal of $B$ contained in the Jacobson radical $\jac(B)$ of $B$. Assume that at least one of the following conditions hold.
\begin{enumerate}[\rm (a)]
\item $\f b$ is contained in $f(A)$.
\item $\f b$ is a torsion $A-$module (with the $A-$module structure inherited by $f$).
\end{enumerate}
Then $\da$ is a total ring of fractions (and it is, in particular, a Pr\"ufer ring).
\end{prop}
\begin{proof}
Let $(a,f(a)+b)$ be a non invertible element of $\da$. The goal is to show that $(a,f(a)+b)$ is a zerodivisor of $\da$. Since $\f b\subseteq \jac(B)$, by Proposition \ref{spec} it follows that
$$
\Max(\da)=\{\f m^{\prime_f}:\f m\in \Max(A)\}.
$$
Thus, there exists a maximal ideal $\f m$ of $A$ such that $(a,f(a)+b)\in\f m^{\prime_f}$, that is $a\in\f m$. Since $A$ is a total ring of  fractions, it follows that $a$ is a zerodivisor of $A$. Hence, we can pick a nonzero element $\alpha\in A$ such that $a\alpha=0$. The following two cases may occur.
\begin{itemize}
\item \emph{Condition (a) holds}. If $\alpha\in {\rm Ann}_A(\f b)$, then it follows immediately that $(a,f(a)+b)(\alpha,f(\alpha))=(0,0)$. Otherwise, let $\beta\in \f b$ be an element such that $f(\alpha)\beta\neq 0$. Since $\f b\subseteq f(A)$, there is an element $x\in f^{-1}(\f b)$ such that $f(x)=\beta$. Of course, $\alpha x\neq 0$ and $(\alpha x,0)\in \da$, since $\alpha x\in f^{-1}(\f b)$. It follows $(a,f(a)+b)(\alpha x,0)=(0,0)$.
\item \emph{Condition (b) holds}. Since $\f b$ is a torsion $A-$module, there exists a regular element $x_0\in A$ such that $f(x_0)b=0$. Of course, $\alpha x_0\neq 0$, since $\alpha\neq 0$. Then $(a,f(a)+b)(\alpha x_0,f(\alpha x_0))=(0,0)$.
\end{itemize}
The conclusion is now clear.
\end{proof}
}
\begin{prop}
We preserve the notation of Proposition \ref{inizio}. The following statements hold
\begin{enumerate}[\rm (1)]
\item If $\da$ is an arithmetical ring, then $A$ is an arithmetical ring.
\item If $\da$ is a Gauss ring, then $A$ is a Gauss ring.
\end{enumerate}
\end{prop}
\begin{proof}
By \cite[Remark 4.6]{da-fi-fo1}, $A$ is a ring retract of $\da$, via the projection $p_{_A}:\da\longrightarrow A$, ($(a,f(a)+b)\mapsto a$). Then, the conclusion follows by applying \cite[Theorem 2.1(1) and Theorem 2.5]{Ba-Ma-Mo}.
\end{proof}
\begin{prop}\label{b_S=0}
We preserve the notation of Proposition \ref{inizio} and Remark \ref{b-loc}. Assume that  $\f b_{S_{\f m}}=\{0\}$, for each $\f m\in \Max(A)\cap V(f^{-1}(\f b))$. Then, the following statements hold.
\begin{enumerate}[\rm (1)]
\item If $A$ is a locally Pr\"ufer ring and  $B_{\f n}$ is a Pr\"ufer ring, for each $\f n\in \Max(B)\setminus V(\f b)$, then $\da$ is a locally Pr\"ufer ring.
\item If $A$ is a Gauss ring and $B_{\f n}$ is a Gauss ring, for each $\f n\in \Max(B)\w V(\f b)$, then $\da$ is a Gauss ring.
\end{enumerate}
\end{prop}
\begin{proof}
By Proposition \ref{spec}, we have
$$
\Max(\da)=\{\f m^{{\prime_{_{\! f}}}} :  \f m\in\Max(A)\}  \cup \{\overline{\f n}^{_{_{f}}} : \f n \in \Max(B)\w V(\f b)\}.
$$
Keeping in mind Proposition \ref{spec}(7) and  that $\f b_{S_{\f m}}=\{0\}$, for each $\f m\in \Max(A)\cap V(f^{-1}(\f b))$, we have that $(\da)_{ \overline{\f n}^{_{_{f}}}}\cong B_{\f n}$, for each $\f n\in \Max(B)\setminus V(\f b)$, and $(\da)_{\f m^{{\prime_{_{\! f}}}}}\cong A_{\f m}$, for each $\f m\in \Max(A)$. Then, statement (1) follows by definition. Statement (2) follows by noting that the property of being Gauss, for a ring, is local.
\end{proof}
\begin{prop}\label{valfib}
We preserve the notation of Definition \ref{fibproddef}, and set $D:=\rho\times_C\sigma$. Let $p_A:D\longrightarrow A$ (resp. $p_B:D\longrightarrow B$) be the restriction to $D$ of the projection of $A\times B$ into $A$ (resp. $B$). The following conditions are equivalent.
\begin{enumerate}[\rm (i)]
\item The set of all ideals of $D:=\rho\times_C\sigma$ is totally ordered by inclusion.
\item At least one of the following statements holds:
\begin{enumerate}[\rm (a)]
\item $\rho$ is injective and the set of all ideals of $p_B(D)$ is totally ordered by inclusion.
\item $\sigma$ is injective and the set of all ideals of $p_A(D)$ is totally ordered by inclusion.
\end{enumerate}
\end{enumerate}
\end{prop}
\begin{proof}
It is immediate that $\Ker(p_A)=\{0\}\times \Ker(\sigma)$ and $\Ker(p_B)=\Ker(\rho)\times \{0\}$.

(ii)$\Longrightarrow$(i). It sufficies to note that, if statement (a) (resp. (b)) holds, then $p_B$ (resp., $p_A$) is an isomorphism of $D$ onto $p_B(D)$ (resp., $p_A(D)$).

(i)$\Longrightarrow$(ii). If the set of all ideals of $D$ is totally ordered by inclusion, obviously each homomorphic image of $D$ has the same property. Thus, if statement (a) is false, $\rho$ is not injective. This implies  $\Ker(p_B)\nsubseteq \Ker(p_A)$, and then we have $\Ker(p_A)\subseteq \Ker(p_B)$, by assumption. It follows immediately that $\sigma$ is injective and $p_A$ is an isomorphism of $D$ onto $p_A(D)$. Thus statement (b) is true.
\end{proof}
\begin{prop}\label{arithm-sur}
We preserve the notation of Proposition  \ref{inizio} and Remark \ref{b-loc}. Assume that for each $\f m\in \Max(A)\cap V(f^{-1}(\f b))$, either the map $f_{\f m}:A_{\f m}\longrightarrow B_{S_{\f m}}$ is surjective or  $f^{-1}(\f b)A_{\f m}\neq \{0\}$.
  Then, the following conditions are equivalent.
\begin{enumerate}[\rm (i)]
\item $\da$ is an arithmetical ring.
\item $A$ is an arithmetical ring, $\f b_{S_{\f m}}=\{0\}$, for each $\f m\in \Max(A)\cap  V(f^{-1}(\f b))$, and, for any $\f n\in \Max(B)\w V(\f b)$, the set of all the ideals of $B_{\f n}$ is totally ordered by inclusion.
\end{enumerate}
\end{prop}
\begin{proof}
(i)$\Longrightarrow$(ii). By \cite[Theorem 1]{Je}, the set of all ideals of each localization of $\da$ at its maximal ideals is totally ordered by inclusion. Thus, Proposition \ref{spec}(7)
implies that in each localization $B_{\f n}$ ($\f n\in\Max(B)\w V(\f b)$) the set of all  ideals is totally ordered by inclusion. Now, let $\f m$ be a maximal ideal of $A$ containing $f^{-1}(\f b)$. By Proposition \ref{spec}(7), the localization $(\da)_{\f m^{\prime f}}$ is isomorphic to $A_{\f m}\Join^{f_{\f m}}\f b_{S_{\f m}}$. If $\pi_{\f m}:B_{S_{\f m}}\longrightarrow B_{S_{\f m}}/\f b_{S_{\f m}}$ is the canonical projection and $\check{f_{\f m}}:=\pi_{\f m}\circ f_{\f m}$, by Proposition \ref{fibprodfarf} the ring $A_{\f m}\Join^{f_{\f m}}\f b_{S_{\f m}}$ is the fiber product of the ring homomorphisms $\check{f_{\f m}}$ and $\pi_{\f m}$. Keeping in mind that ${f_{\f m}}^{-1}(\f b_{S_{\f m}})=f^{-1}(\f b)A_{\f m}$ (Remark \ref{b-loc}) and applying Proposition \ref{valfib}, it follows that $\f b_{S_{\f m}}=\{0\}$. Thus, by Proposition \ref{spec}(7), $A_{\f m}$ is isomorphic to $(\da)_{\f m^{\prime f}}$, for each maximal ideal $\f m$ of $A$. This proves that $A$ is an arithmetical ring.

(ii)$\Longrightarrow$(i). Apply \cite[Theorem 1]{Je}, Proposition \ref{inizio}(3) and the local structure of $\da$ (Proposition \ref{spec}(7)).
\end{proof}
\begin{prop}\label{wgldim-sur}
We preserve the notation of Proposition \ref{inizio} and Remark \ref{b-loc}. Assume that, for each maximal ideal $\f m$ of $A$ containing $f^{-1}(\f b)$, either the map $f_{\f m}:A_{\f m}\longrightarrow B_{S_{\f m}}$ is surjective or $f^{-1}(\f b)A_{\f m}\neq \{0\}$. Then, the following conditions are equivalent.
\begin{enumerate}[\rm (i)]
\item $\da$ has weak global dimension at most 1.
\item $A$ has weak global dimension at most 1, $B_{\f n}$ is a valuation domain, for each  $\f n\in\Max(B)\w V(\f b)$ and $\f b_{S_{\f m}}=\{0\}$, for each $\f m\in \Max(A)\cap V(f^{-1}(\f b))$.
\end{enumerate}
\end{prop}
\begin{proof}
(i)$\Longrightarrow$(ii). By Proposition \ref{spec}(7), we have $(\da)_{\overline{\f n}^f}\cong B_{\f n}$, for any maximal ideal $\f n$ of $B$ not containing $\f b$. Then, it follows, by definition, that $B_{\f n}$ is a valuation domain for each $\f n\in \Max(B)\w V(\f b)$.
Now, let $\f m$ be a maximal ideal of $A$ containing $f^{-1}(\f b)$. Since, in particular, $\da$ is an arithmetical ring, it follows  $\f b_{S_{\f m}}=\{0\}$, by Proposition \ref{arithm-sur}. Thus, by  Proposition \ref{inizio}(3), the localization $A_{\f m}$ is isomorphic to the the valuation domain  $(\da)_{\f m^{\prime_f}}$, for any maximal ideal $\f m$ of $A$. This proves that  $A$ has weak global dimension $\leq$ 1.

(ii)$\Longrightarrow$(i). Apply the local structure of $\da$ (Proposition \ref{spec}(7)). Note that (ii) implies (i), without any extra assumption.
\end{proof}
%

To give conditions to make $\da$ a semi--hereditary ring, we want to use the following characterization.
\begin{thm}\label{semih-carat}
{\rm (\cite[Corollary 4.2.19]{Glaz})}
Let $A$ be a ring. Then, $A$ is semi--hereditary if and only if $A$ is coherent and the weak global dimension of $A$ is at most 1.
\end{thm}
Let $\phi:A\longrightarrow B$ be a ring homomorphism and let $M$ be a $B-$module. We shall denote by $\cdot_\phi$ the scalar multiplication, induced by $\phi$, making $M$ an $A-$module.
{\bc
\begin{lem}\label{finite}
We preserve the notation of Proposition \ref{inizio}. If $\f b$ is a finitely generated $A-$module (with the $A-$module structure induced by $f$), then the ring embedding $\iota: A\longrightarrow \da$ is finite.
\end{lem}
\begin{proof}
Let $\{b_1,\z,b_n\}\subseteq \f b$ be a finite set of generators of the $A-$module $\f b$, and fix an element $(a,f(a)+b)\in \da$. Then, there exist elements $a_1,\z,a_n\in A$ such that $b=\sum_{i=1}^na_i\cdot_fb_i=\sum_{i=1}^nf(a_i)b_i$. It follows immediately that
$$
(a,f(a)+b)=a\cdot_\iota(1,1)+\sum_{i=1}^n
a_i\cdot_\iota
(0,b_i).
$$
This proves that $\{(1,1),(0,b_1),\z,(0,b_n)\}\subseteq \da$ is a finite set of generators of $\da$ as an $A-$module (with the structure induced by $\iota$), i.e. $\iota$ is finite.
\end{proof}
}
\begin{prop}\label{cohe-farf}
We preserve the notation of Proposition \ref{inizio}. Then, the following statements hold.
\begin{enumerate}[\rm (1)]
\item If $\da$ is a coherent ring, then $A$ is coherent.
\item If $A$ is  a coherent ring and $\f b$ is a coherent $A-$module (with the structure induced by $f$), then  $\da$ is a coherent ring.
\end{enumerate}
\end{prop}
\begin{proof}
Statement (1) follows by \cite[Remark 4.6]{da-fi-fo1} and \cite[Theorem 4.1.5]{Glaz}.

(2). We begin by noticed that, since $\f b$ is, in particular, a finitely generated $A-$module, the ring embedding $\iota$ is finite, by Lemma \ref{finite}. Now, let $p_{_A}:\da\longrightarrow A,\,\,p_{_B}:\da \longrightarrow B$ be the projections. Then, $p_{_A}$ (resp. $p_{_B}$) induces on $A$ (resp. $\f b$) a structure of $\da-$module. With these structures, we have the following short exact sequence
$$0\longrightarrow \f b\stackrel{i }{\longrightarrow}\da\stackrel{p_{_A}}
{\longrightarrow}A
\longrightarrow 0,$$
of $\da-$modules, where $i:\f b\longrightarrow \da$ is  defined by $\beta\mapsto (0,\beta)$, for each $\beta\in \f b$. Let $\iota:A \hookrightarrow \da$ be the ring enbedding such that $a\mapsto (a,f(a))$, for each $a\in A$. On the $\da-$module $\f b$, the map $\iota$ induces the following scalar multiplication
$$
a\cdot_\iota \beta:=(a,f(a))\cdot_{p_{_B}}\beta=
p_{_B}((a,f(a)))
\beta=f(a)\beta \qquad (a\in A,\beta\in \f b)
$$
It follows that the structure of $A-$module given to $\f b$ by $\iota$ is the same structure induced on $\da$ by $f$. Since $\iota$ is finite and  $\f b$ is a coherent $A-$module, by \cite[Corollary 1.1]{har} it follows that $\f b$ is a coherent $\da-$module. Moreover, $\iota$ induces to the $\da-$module $A$ the following scalar multiplication
$$
a\cdot_\iota\alpha:=(a,f(a))\cdot_{p_{_A}}
\alpha=p_{_A}((a,f(a)))\alpha= a\alpha \qquad (a,\alpha\in A)
$$
Thus $\iota$ induces on $A$ its natural structure of module over itself. Since $A$, by assumption, is a coherent ring, it follows that it is a coherent $\da-$module, again by \cite[Corollary 1.1]{har}. Then $\da$ is a coherent $\da-$module, by \cite[Pag. 43, Exercise 11(a)]{Bour}, that is, $\da$ is a coherent ring.
\end{proof}
{\bc
\begin{cor}\label{semihsufficient}
We preserve the notation of Proposition \ref{inizio} and Remark \ref{b-loc}, and assume that $\f b_{S_{\f m}}=\{0\}$, for each maximal ideal  $\f m$ of $A$ containing $f^{-1}(\f b)$. If $A$ is a semi--hereditary ring (resp. semi--hereditary and Noetherian ring), $B_{\f n}$ is a valuation domain, for each $\f n\in\Max(B)\w V(\f b)$ and $\f b$ is a coherent $A-$module (resp. finitely generated $A-$module), with the structure induced by $f$, then $\da$ is a semi--hereditary ring.
\end{cor}
\begin{proof}
Apply Theorem \ref{semih-carat}, Proposition \ref{wgldim-sur}((ii)$\Longrightarrow$(i)) and Proposition \ref{cohe-farf}, keeping in mind that, if $A$ is a Noetherian ring, an $A-$module is coherent if and only if it is finitely generated.
\end{proof}
Recall that an integral domain $A$ is \emph{almost Dedekind} if $A_{\f m}$ is a DVR for each maximal ideal $\f m$ of $A$. Thus, in particular, an almost Dedekind domain is a Pr\"ufer domain.
\begin{ex}
Let $A$ be a non-Noetherian almost Dedekind domain having at least two distinct principal maximal ideals $\f m:=(m),\f n:=(n)$ (such a domain exists, see \cite{fo-hu-pa}), set $B:=A/(\f m \cap \f n)$, let $f:A\longrightarrow B$ be the canonical projection and set $\f b:=\f m/(\f m\cap \f n)$. Trivially, $f^{-1}(\f b)=\f m$ and, since $f(n)\in S_{\f m}$, it follows that $\f b_{S_{\f m}}=\{0\}$. Let $\overline{\f n}:=\f n/(\f m\cap \f n)$ be the unique maximal ideal of $B$ not containing $\f b$. Obviously, the localization $B_{\overline{\f n}}$ is isomorphic to the field $A/\f n$. Moreover, the natural map $p:A\longrightarrow \f b$, $a\mapsto f(am)$, is clearly $A-$linear,  surjective  and $\Ker(p)=\f n$. This shows that $\f b$ is finitely presented as an $A-$module. Then, keeping in mind that $A$ is a coherent ring, being it a Pr\"ufer domain, and applying  \cite[Exercise 12 (a)($\beta$)]{Bour}, it follows that $\f b$ is a coherent $A-$module. Then $\da$ is a semi--hereditary ring, by Corollary \ref{semihsufficient}.
\end{ex}
\begin{ex}
Preserve the notation of Proposition \ref{inizio}. The fact that $\da$ is semi--hereditary does not imply, in general, that $\f b$ is coherent as an $A-$module and $\f b_{S_{\f m}}=\{0\}$, for each $\f m\in\Max(A)\cap V(f^{-1}(\f b))$. For example, let $T$ be an indeterminate over $\mathbf Q$, and let $A:=\mathbf Z$, $B:=\mathbf Q[T]$, $\f b:=T\mathbf Q[T]$, $f:A\longrightarrow B$ be the inclusion. Then $\da$ is isomorphic to the ring $\mathbf Z+T\mathbf Q[T]$, by \cite[Example 2.5]{da-fi-fo1}. Moreover, by \cite[Theorem 1.3]{Ho-Ta}, it follows easily that $\da$ is a Pr\"ufer domain (i.e. a semi--hereditary domain). But, clearly, $\f b$ is not finitely generated as an $A-$module and $\f b_{S_{\f m}}\neq \{0\}$, for each $\f m\in\Max(A)$.
\end{ex}
\begin{cor}\label{semih-sur}
We preserve the notation of Proposition \ref{inizio} and Remark \ref{b-loc}. Assume that $\f b$ is a coherent $A-$module and that, for each $\f m\in\Max(A)\cap V(f^{-1}(\f b))$, either $f_{\f m}$ is a surjective ring homomorphism or $f^{-1}(\f b)A_{\f m}\neq\{0\}$. Then, the following conditions are equivalent.
\begin{enumerate}[\rm (i)]
\item $\da$ is a semi--hereditary ring.
\item $A$ is a semi--hereditary ring, $B_{\f n}$ is a valuation domain, for each $\f n\in \Max(B)\w V(\f b)$ and $\f b_{S_{\f m}}=\{0\}$, for each $\f m\in \Max(A)\cap V(f^{-1}(\f b))$.
\end{enumerate}
\end{cor}
\begin{proof}
(ii)$\Longrightarrow$(i). It is the statement of Corollary \ref{semihsufficient}.

(i)$\Longrightarrow$(ii). By Proposition \ref{cohe-farf}(1), $A$ is a coherent ring. Then, it sufficies to apply Theorem \ref{semih-carat} and Proposition \ref{wgldim-sur} to complete the proof.
\end{proof}
}


\end{document}